\documentclass{amsart}
\usepackage{hyperref}
\usepackage[numeric,short-months]{amsrefs}
\usepackage[all]{xy}
\usepackage{math}

\newcommand\dunderline[1]{\underline{\underline{#1}}}

\newtheorem{thmp}[thm]{\dunderline{Theorem}}
\newtheorem{corp}[thm]{\dunderline{Corollary}}

\begin {document}
\author{Laurent Bartholdi}
\author{Illya I. Reznykov}
\date{4 January 2007; compiled \today}
\title{A Mealy machine with polynomial growth of irrational degree}

\begin{abstract}
  We consider a very simple Mealy machine (two non-trivial states over
  a two-symbol alphabet), and derive some properties of the semigroup
  it generates. It is an infinite, finitely generated semigroup, and
  we show that the growth function of its balls behaves asymptotically
  like $n^\alpha$, for $\alpha = 1+{\log 2} / {\log
    \frac{1+\sqrt{5}}{2}}$; that the semigroup satisfies the identity
  $g^6=g^4$; and that its lattice of two-sided ideals is a chain.
\end{abstract}

\maketitle


\section{Introduction}\label{sec:introduction}
Algebraic objects may be defined by universal properties, like the
identities they satisfy, or a presentation as a quotient of two free
objects. They may also be defined by their action on simpler objects,
viz, for an algebra, as endomorphisms of a vector space; for a group,
as permutations of a set, and for a semigroup, as mappings of a set.

If the semigroup to be defined is infinite, it naturally must act on
an infinite set; extra conditions must be imposed on the action to
ensure that it remains describable (by a finite sequence of
mathematical symbols, say). A specially interesting class of
semigroups appear by requesting that the infinite set be the set $X^*$
of words over a finite alphabet, and that the action be given by a
finite-state automaton.

The growth function $\gamma(\ell)$ of a finitely generated semigroup
$\Gamma$ --- the number of semigroup elements that can be obtained as
products of at most $\ell$ generators --- is an important invariant of
the semigroup.  It depends on the chosen generating set, but its
asymptotics do not.  This function is at most exponential. If it
is bounded by a polynomial, then $\Gamma$ is of \emph{polynomial
  growth}; following~\cite{gelfand-k:dimension}, its
\emph{Gelfand-Kirillov dimension} (abbreviated GK dimension) is then
defined as the infimum of those $\alpha$ such that
$\gamma(\ell)/\ell^\alpha$ does not converge to $0$.

The ``Bergman gap''~\cite{krause-l:gkdim}*{Theorem~2.5} asserts that a
semigroup may have GK dimension $0$, $1$, or $\ge2$; by a result by
Warfield~\cite{warfield:tensor}, there exist semigroups of GK
dimension $\beta$ for any $\beta\ge2$.  Belov and
Ivanov~\cites{belov-i:growth1,belov-i:growth2} construct finitely
presented semigroups with non-integer GK dimension.
Shneerson~\cite{shneerson:interm} constructs relatively free (i.e.,
free relative to an identity) semigroups of intermediate growth
(asymptotically $m^{\log{m}}$). His proof involves Fibonacci numbers,
as does the present construction.

\begin{figure}
  \[\xymatrix{
    *+[o][F]{f}\ar[rr]^{0\to0} \ar@(dl,ul)[]^{1\to0} &&
    *+[o][F]{s}\ar@/^1pc/[rr]^{0\to1} \ar@/_1pc/[rr]_{1\to0} &&
    *+[o][F]{e}\ar@(u,r)[]^{0\to0} \ar@(d,r)[]_{1\to1}
  }\]
  \caption{The automaton $I$}\label{fig:automaton}
\end{figure}
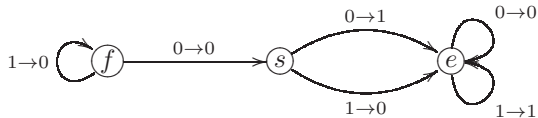

This paper describes the semigroup $\Gamma(I)$ generated by the
$3$-state automaton $I$ in Figure~\ref{fig:automaton}
(see~\S\ref{sec:definitions} for the definition of a semigroup
generated by an automaton), and proves that its GK dimension is
irrational. A much more precise statement appears in
Theorem~\ref{thm:growth}.

The main purpose of this paper is to show that a very simple-minded
finite-state automaton can produce a semigroup with a highly unusual
growth pattern, asymptotically $n^\alpha$ for an irrational $\alpha$.
Other exotic types of growth of automata have also been discovered;
for example, of type $e^{\sqrt n}$ in~\cite{bartholdi-r-s:interm2x2},
and of type $n^{\log(n)/2\log m}$ for integer $m$,
in~\cite{reznykov-s:3state}. It is intriguing that in all these cases
the determination of the growth function relies on enumerations of
partitions with certain constraints.

\subsection{Plan}
The next section gives an automata-free version of the semigroup
studied in this paper, along with a presentation of its main results.
Section~\ref{sec:definitions} gives the necessary definitions of
semigroups generated by automata. Statements whose proof are too long
to fit smoothly in the text are stated as \dunderline{underlined}.
Their proofs appear in Section~\ref{sec:proofs}.

\subsection{Notation}
All actions are written on the right in this paper. The identity is
written $e$. We use $=$ for equality of group elements, and $\equiv$
for graphical (letter-by-letter) equality of words representing these
elements. We denote by $X^*$ the free monoid on the set $X$. The
length of a word $w$ is denoted by $\|w\|$; this is also the length of
a minimal word representing a semigroup element.

The integers are written $\Z$, and the naturals (containing $0$) are
written $\N$. Congruence is written $\equiv$, and the `mod 2' operation
(remainder after division by $2$) is written $n\%2$ as in the C
programming language.

Sequences of the form $f_af_{a+2}\dots f_b$, sometimes simply written
$f_a\dots f_b$, appear throughout the paper. They are taken to be $e$
if $a>b$, and $f_af_{a+2}f_{a+4}\dots f_{b-2}f_b$ otherwise.

\subsection{Acknowledgments} We are very grateful to the anonymous
referee who suggested substantial improvements to the paper, in
particular a clarification of the role of rewriting systems.

\section{Main results}\label{sec:main}
Consider the following transformations $s,f$ of the integers $\Z$:
\begin{equation}\label{eq:zaction}
  \begin{aligned}
    x^s&=\begin{cases} x-1 & \text{ if $x$ is odd,}\\
      x+1 & \text{ if $x$ is even;}\end{cases}\\
    x^f&=\begin{cases} x-2^n-1 &
      \text{ if there exists $n\ge0$ such that }x\equiv3\cdot2^n\pmod{2^{n+2}},\\
      x+3\cdot2^n-1 &
      \text{ if there exists $n\ge0$ such that }x\equiv2^n\pmod{2^{n+2}},\\
      -1 & \text{ if }x=0.
    \end{cases}
  \end{aligned}
\end{equation}
Note that $f$ is uniquely defined, because every $x\in\Z$ is either
$0$ or of the form $2^ny$ for unique $n\ge0$ and odd $y$, with either
$y\equiv1$ or $y\equiv3$ modulo $4$.

\subsection{The semigroup} Let $F$ be the semigroup generated by the
transformations $s$ and $f$.  Define furthermore the elements $f_n$ of
$F$ by $f_1=s$, $f_2=f$, and inductively $f_n=f_{n-2}f_{n-1}$ for
$n\ge3$. For example, $f_3=sf$, $f_4=fsf$, $f_5=sf^2sf$. These words
are sometimes called the ``Fibonacci sequence'' ---
see~\S\ref{ss:normalform} for the connection.

Recall~\cite{epstein-:wp} that a rewriting system for a semigroup
$\Gamma$ generated by a set $Q$ is a set of equations, called
\emph{rules}, of the form $\ell\to r$, with $\ell,r\in Q^*$. An
elementary reduction in a word $w\in Q^*$ is the replacement of a
subword equal to the left-hand side of a rule by the right-hand side;
if no elementary reduction is possible, the word is \emph{reduced}. A
rewriting system is \emph{terminating} if there is no $w\in Q^*$ to
which an infinite sequence of elementary reductions can be applied;
and it is \emph{confluent} if the reduced word obtained after applying
as many elementary reductions as possible does not depend on the
choice of elementary reductions. It is \emph{complete} if it is
terminating and confluent; the set of reduced words is then in
bijection with the semigroup, through the natural evaluation map
$Q^*\to\Gamma$, and is called a \emph{normal form} for $\Gamma$.

For $n\ge1$, define words $r_n$ and $r'_n$ over the alphabet $\{s,f\}$
as follows:
\begin{align}
  r_1 &=s^2,\qquad r'_1=e,\notag\\
  \label{eq:def:r_n} r_n &= f_{n+1}f_n^2,\\
  \label{eq:def:r'_n} r'_n &= f_{n\%2+1}(f_{n\%2+5}f_{n\%2+7}\cdots f_{n-1}f_{n+1})f_n.
\end{align}
In particular, $r_2=sf^3$ and $r'_2=sf$, and $r_3=f_4f_3^2$ and
$r'_3=f_4$. For $n\ge3$, we have $\|r'_n\|=\|r\|-4$.

\begin{thmp}\label{thm:pres}
  The semigroup $F$ is infinite, and admits as presentation
  \[\langle s,f|\,r_n = r'_n\text{ for all }n\ge1\rangle,\]
  Furthermore, after the relations $sw=sw'$ (which occur for even $n$)
  are replaced by the equivalent relations $w=w'$, this presentation
  is a complete rewriting system.
\end{thmp}

\subsection{A normal form}\label{ss:normalform}
Let $\Phi_n$ denote the sequence of Fibonacci numbers, defined by
$\Phi_1=\Phi_2=1$ and $\Phi_n=\Phi_{n-2}+\Phi_{n-1}$ for $n\ge3$. Then
\begin{thmp}\label{thm:nf}
  Every element $g$ of $F$ admits a unique representation as a word of
  the form
  \begin{equation}\label{eq:nf}
    w_g=s^\epsilon f_{i_1}f_{i_2}\cdots f_{i_m}\cdots f_{i_n},
  \end{equation}
  for some $n\ge0$, some $\epsilon\in\{0,1\}$, and some indices
  $i_1\dots,i_n$ satisfying
  \[3\le i_1,i_1+1<i_2,\dots,i_{m-1}+1<i_m>i_{m+1}>\cdots>i_n\ge1.\]
  We call $i_m$ the \emph{maximal index} of $w_g$.

  When spelled out in the generating set $\{s,f\}$, this
  representation of $g$ is essentially minimal: if $i_1$ is even, then
  $w_g$ is the unique minimal-length representation of $g$, while if
  $i_1$ is odd, then an initial $s^2$ must be cancelled to obtain the
  unique minimal representation of $g$. Its length
  (see~\S\ref{ss:metrics}) is
  \[\|g\|=(-1)^{i_1}\epsilon+\Phi_{i_1}+\cdots+\Phi_{i_n}.\]
\end{thmp}

By a slight extension of the definition of rewriting system, let us
admit rules of the form $\wedge\ell\to r$ that mean that the subword
$\ell$ may be replaced in $w$ by the word $r$ if it is a prefix of
$w$. We will actually consider the following rewriting system; it is
on an infinite generating set, but its rules are much simpler, since
their left-hand sides have length at most $3$.
\begin{thmp}\label{thm:rws}
  On the generating set $\{f_i:i\ge1\}$, the semigroup $F$ admits a
  complete rewriting system with rules
  \begin{align}
    f_1^2&\longrightarrow e\tag{N1}\\
    f_af_{a+1}&\longrightarrow f_{a+2}\quad(a\ge1)\tag{N2}\\
    f_a^2&\longrightarrow f_{a-2}f_{a+1}\quad(a\ge3)\tag{N3}\\
    f_af_2^2&\longrightarrow f_a\quad(a\ge2)\tag{N4}\\
    \wedge f_2&\longrightarrow f_1f_3\tag{N5}\\
    f_{a+1}f_af_{a+3}&\longrightarrow f_{a+3}f_{a+2}\quad(a\ge1)\tag{N6}\\
    f_{a+2}f_af_{a+3}&\longrightarrow f_af_{a+3}f_{a+2}\quad(a\ge1)\tag{N7}\\
    \tag{N8} f_{a+p}f_af_{a+q}&\longrightarrow \underbrace{f_{a+p-2}\cdots
      f_a}_{p-1\text{ terms}}f_2^2\cdots
    f_{a-2}^2f_{a-1}f_{a+q}\\
    &(a\ge2,p\ge1,q\ge2\text{ even})\notag\\
    f_{a+p}f_af_{a+q}&\longrightarrow \underbrace{f_{a+p-2}\cdots
      f_{a+2}}_{p-3\text{ terms}}f_2^2\cdots
    f_a^2f_{a+q}\quad(a\ge1,p\ge3,q\ge3\text{ odd})\tag{N9}\\
    f_{1+p}f_1f_{1+q}&\longrightarrow f_{1+p}f_2f_4\cdots
    f_q\quad(p\ge1,q\ge2\text{ even})\tag{N10}.
  \end{align}
\end{thmp}

From now on, we identify $g$ with $w_g$, and call the word $w_g$ the
\emph{normal form} of $g$. We will never use a different notation for
a semigroup element and its normal form.

\subsection{Growth} Using this minimal representation, the (ball)
growth function of $F$,
\[\gamma(\ell)=\#\{g\in F:\,\|g\|\le\ell\},\]
may be quite precisely estimated; namely
\begin{thmp}\label{thm:growth}
  There are constants $C,D>0$ such that the growth function of $F$
  satisfies
  \[C\ell^\alpha\le\gamma(\ell)\le D\ell^\alpha\] for all $\ell\in\N$,
  where $\varphi=(1+\sqrt5)/2$, and
  $\alpha=1+\log2/\log\varphi\approx2.4401$, and
  $C=\frac{2\sqrt5^\alpha}{\sqrt
    5\varphi^2(2\varphi-1)(2\varphi)^\alpha}$, and
  $D=\frac{2\sqrt5^\alpha}{\sqrt 5\varphi^2(2\varphi-1)}$.

  Therefore, the Gelfand-Kirillov dimension of $F$ is $\alpha$.
\end{thmp}
Experimental computations indicate that actually the function
$\gamma(\ell)/\ell^\alpha$ does not converge, but oscillates between
$\approx0.201$ and $\approx0.205$, reaching its maxima at Fibonacci
numbers.

\subsection{Identities} We have not determined the complete set of
identities satisfied by $F$; nevertheless, we know
\begin{thmp}\label{thm:identity}
  The semigroup $F$ satisfies the identity $g^6=g^4$.
\end{thmp}

\subsection{The ideal structure of $F$}
We describe in this subsection the quotient and ideal structure of
$F$. First, we may consider for all $n\in\N$ the quotient of $F$,
denoted $W_n$, acting as transformations on $X^n$. Let us denote, for
$n\in\N$,
\[z_n = \begin{cases}
  f_3 f_5 \cdots f_{n+2} & \text{ if $n$ is odd,} \\
  f_4 f_6 \cdots f_{n+2} & \text{ if $n$ is even.}
\end{cases}
\]
Then we have
\begin{thmp}\label{thm:quotients}
  \begin{enumerate}
  \item The semigroup $W_n$ is presented as follows:
    \[
    W_n = \left\langle s,f \, \vline \, r_k = r'_k, 1 \le k \le n+2,
      \text{ and } sz_n = fz_n = z_n, f_{n+2} f_{n+1} = z_n
    \right\rangle.
    \]
  \item The elements of $W_n$ may be described as all normal forms
    $f_1\cdots f_{i_m}\cdots f_{i_t}$ whose maximal index $i_m$ is $< n
  +2$, normal forms of the form $z_n f_{i_{m-1}} \cdots f_{i_t}$
    for $n \ge i_{m-1} > \cdots > i_t$, and those normal forms of
    maximal index $n+2$ that include neither $z_n$ nor $f_{n+2}
    f_{n+1}$.
  \item The semigroup $W_n$ has order $\sum_{k = 1}^n 2^{k+2}
    \Phi_{k}-2^n-2^{n+1} \Phi_{n}+2$.
  \end{enumerate}
\end{thmp}
We identify for the remainder of this section the quotient semigroup
$W_n$ with the set of normal forms of its elements, which then form a
subset of $F$.

\begin{cor}
  The Hausdorff dimension (in the sense of~\S\ref{ss:hausdorff}) of
  $F$ is $0$.
\end{cor}
\begin{proof}
  Using the count of elements in $W_n$, we have
  \[
  \#W_n \le 2+\sum_{k = 1}^n 2^{k+2} \Phi_{k}
  \le\frac{4}{\sqrt 5}\frac{(2\varphi)^{n+1}-1}{2\varphi-1}\le(2\varphi)^{n+1},
  \]
  so
  \[
  \operatorname{Hdim}(F) = \liminf_{n\to\infty}
  \frac{\log\#W_n}{\log\#\operatorname{End}(X^n)} =
  \liminf_{n\to\infty}
  \frac{\log(2\varphi)^{n+1}}{\log2^{2(2^n-1)}} = 0.\qedhere
  \]
\end{proof}

\begin{corp}\label{cor:trace}
  Let $g$ be an arbitrary element of $F$ with maximal index $n \ge 3$.
  Then its trace (in the sense of~\S\ref{ss:traces}) satisfies
  \[\tau(g) = \begin{cases} 2^{2-n} & \text{ if $g$ includes $z_{n-2}$
      or $f_n f_{n-1}$ in its normal form,} \\
    2^{3-n} & \text{ otherwise.}
  \end{cases}
  \]
  The ideal $I_n=Ff_{n+3}F$ coincides with $F_{2^{-n}}$, as defined
  in~\S\ref{ss:traces}.
\end{corp}

\begin{thmp}\label{thm:ideals}
  All two-sided ideals of $F$ are of the form $Ff_nF$, and
  equivalently of the form $F_\xi$ as defined in~\S\ref{ss:traces}.
  They therefore form a unique chain.
\end{thmp}

Recall that if $I$ is an ideal of $F$, the quotient $F/I$ is the
semigroup whose elements are the equivalence classes of the congruence
$I\times I\cup\Delta$, where $\Delta$ is the identity relation (=
diagonal) of $F$.
\begin{prop}
  For $n\ge2$, the quotient $F/I_n$ has order $\sum_{k=1}^n 2^{k+2}
  \Phi_k-2^{n+1}-2^{n+1}\Phi_n+3$. It is obtained from $W_n$ by
  identifying together all idempotents of rank $1$, i.e.\ all maps
  $f_p:X^n\to X^n$ defined by $f_p(x_1\dots x_n)=p$, for all $p\in
  X^n$.
\end{prop}
\begin{proof}
  By Theorem~\ref{thm:quotients}, the quotient $F/I_n$ is a quotient
  of $W_n$. The transformation $f_{n+3}$ induces a transformation of
  rank $1$ on $X^n$, and its left multiples, which are identified in
  $F/I_n$, are the maps $f_p$ in the statement of the Proposition.
\end{proof}

\section{Definitions}\label{sec:definitions}
A \emph{finite-state automaton}, also called \emph{Mealy machine}, is
comprised of the following data: a set $X$, the \emph{alphabet}; a set
$Q$, the \emph{states}; a function $\tau:X\times Q\to Q$, the
\emph{transition}; a function $\pi:X\times Q\to X$, the \emph{output}.
In this section we suppose that the alphabet is $X=\{0,\dots,d\}$; we
will later specialize to $d=1$.

Such an automaton $\mathcal A$ is usually represented, as in the
figure above, as a graph. The states are represented by vertices, and
there is an edge from vertex $q$ to vertex $r$, labelled $i\to o$,
whenever $\tau(i,q)=r$ and $\pi(i,q)=o$.

\subsection{Action} States $q\in Q$ of $\mathcal A$ yield
transformations of $X^*$, which are defined simultaneously as follows:
\begin{equation}\label{eq:action}
  (i_1\dots i_n)^q=\pi(i_1,q)(i_2\dots i_n)^{\tau(i_1,q)}.
\end{equation}
These transformations generate a semigroup $\Gamma(\mathcal A)=\langle
Q\rangle$, called the \emph{semigroup generated by $\mathcal A$.}

\subsection{Decomposition} Yet another way of describing the automaton
$\mathcal A$ is via its \emph{semigroup decomposition}.  Given
$g\in\Gamma$, it yields by~\eqref{eq:action} a transformation $\pi_g$
of $X$ by restriction to length-$1$ words; and, for all $i\in X$, a
transformation $iX^*\to i^{\pi_g}X^*$, again by restriction.
By~\eqref{eq:action}, the composition $X^*\to iX^*\to i^{\pi_g}X^*\to
X^*$ is again an element of $\Gamma$, which we denote by $g_i$. We
write
\begin{equation}\label{eq:gendecomposition}
  \phi(g)=\pair<g_0,\dots,g_d>\pi_g
\end{equation}
for the decomposition of $g\in\Gamma$, with $g_i\in\Gamma$ and
$\pi_g:X\to X$. Multiplication of such decompositions obeys the rule
\begin{equation}\label{eq:product}
  \pair<g_0,\dots,g_d>\pi_g\pair<h_0,\dots,h_d>\pi_h
  =\pair<g_0h_{\pi_g(0)},\dots,g_dh_{\pi_g(d)}>\pi_g\pi_h.
\end{equation}
It is therefore sufficient to know the decomposition of generators,
and these are determined by the transition and output functions
$\tau,\pi:X\times Q\to Q$ of $\mathcal A$, by
\[\phi(q)=\pair<\tau(0,q),\dots,\tau(d,q)>(i\mapsto\pi(i,q)\,\forall i).\]

\subsection{Metrics}\label{ss:metrics} On the semigroup $\Gamma$
generated by a set $Q$, define the ``norm''
\[\|g\|=\min\{\ell\in\N:\,g=q_1\cdots q_\ell,\quad q_i\in Q\,\forall i\}.\]
(we have $\|gh\|\le\|g\|+\|h\|$, which justifies calling it a norm.)
The \emph{ball growth function} is the function $\gamma:\N\to\N$
defined by
\[\gamma(\ell)=\#\{g\in\Gamma:\,\|g\|\le \ell\};\]
it measures the volume growth of balls in the discrete normed space
$\Gamma$.

\subsection{Contraction} An automaton is
\emph{contracting}~\cite{grigorchuk:burnside} if there are constants
$C$ and $\eta<1$ such that whenever
$\phi(g)=\pair<g_0,\dots,g_d>\pi_g$, we have
\[\|g_i\|\le\eta\|g\|+C.\]
Contracting automata are most studied, in part because contraction
gives a natural strategy of proof by induction.

\subsection{Traces}\label{ss:traces} Let $g$ be a transformation of
$X^*$ given by an automaton. Then for all $n\in\N$ we have
$\#g(X^{n+1})\le\#X\#g(X^n)$, so the limit
\[\tau(g) := \lim_{n\to\infty}\frac{\#g(X^n)}{\#X^n}\in[0,1]\]
exists. We call it the \emph{trace} of $g$. The following facts are
easily checked:
\begin{itemize}
\item $\tau(gh)\le\min\{\tau(g),\tau(h)\}$;
\item $\tau(g)=1$ if and only if $g$ is invertible.
\end{itemize}
Therefore, if $\Gamma$ is the semigroup generated by an automaton, we
may define for any $\xi\in[0,1]$ a two-sided ideal
$\Gamma_\xi=\{g\in\Gamma:\,\tau(g)\le\xi\}$. We have
$\Gamma_\xi=\Gamma$ if and only if $\xi=1$.

\subsection{Hausdorff dimension}\label{ss:hausdorff} Let $\Gamma$
be a semigroup acting on a rooted tree $X^*$. There is another metric
on $\Gamma$, defined by
\[d(g,h)=\exp(-\max\{n:\,v^g=v^h\text{ for all }v\in X^n\}),\] with
the convention that $\exp(-\infty)=0$. This turns $\Gamma$ into a
metric space of diameter at most $1$. The \emph{Hausdorff dimension}
of $\Gamma$ is defined as the Hausdorff dimension of this metric
space. Let $\operatorname{End}(X^n)$ denote the set of
prefix-preserving maps $X^n\to X^n$. If $\#X=d+1$, then
\[\#\operatorname{End}(X^n)=((d+1)^{d+1})^{((d+1)^n-1)/d}.\]
Then the Hausdorff dimension of $\Gamma$ is given by the formula
\[\operatorname{Hdim}(\Gamma)=\liminf \limits_{n \to \infty}
\frac{\log\#\Gamma_n}{\log\#\operatorname{End}(X^n)},
\]
where $\Gamma_n$ is the quotient of $\Gamma$ acting as a
transformation semigroup of $X^n$; see~\cite{abercrombie:subgroups}
for an analogous definition in the case of profinite groups.

\section{The automaton $I$ and the semigroup $F$}
We now consider the automaton $I$ from Figure~\ref{fig:automaton}.
This automaton has three states $f,s,e$ and a two-letter alphabet
$X=\{0,1\}$.

Consider the following transformations $\sigma,\zeta$ of $X$:
\[\sigma(i)=1-i,\qquad \zeta(i)=0.\]
Then the decomposition of the states are
\[\phi(f)=\pair<s,f>\zeta,\quad\phi(s)=\pair<e,e>\sigma,\quad
\phi(e)=\pair<e,e>.\]
Clearly $e$ acts as the identity transformation, and $s$ is
invertible, of order $2$.

\subsection{Action on integers} We show now that the semigroup
defined in Section~\ref{sec:main} by its action on integers is the
semigroup generated by an automaton.
\begin{thm}
  The semigroups $F$ and $\Gamma(I)$ are isomorphic.
\end{thm}
\begin{proof}
  The action of $F=\langle s,f\rangle$ given
  in~\eqref{eq:zaction} extends to a continuous
  action of $F$ on the $2$-adics $\Z_2$. Consider the following
  bijection $\Theta$ between $X^\infty$ and $\Z_2$:
  \[\Theta(x_1x_2\dots)=\sum_{i=1}^\infty (1-x_i)2^{i-1}.\]

  Let us denote temporarily $\Gamma(I)=\langle\tilde s,\tilde
  f\rangle$. The theorem will follow from $\Theta(x^{\tilde
    s})=\Theta(x)^s$ and $\Theta(x^{\tilde f})=\Theta(x)^f$ for all
  $x\in X^\infty$.

  Consider therefore $x=x_1x_2\dots\in X^\infty$, and write
  $y=x_2x_3\dots$.  If $x_1=0$, then
  \[\Theta(x)^s=\Theta(0y)^s=(1+2\Theta(y))^s=2\Theta(y)=\Theta(1y)=\Theta(x^{\tilde s}),\]
  while if $x_1=1$, then
  \[\Theta(x)^s=\Theta(1y)^s=(2\Theta(y))^s=1+2\Theta(y)=\Theta(0y)=\Theta(x^{\tilde s}).\]

  Let next $n\in\N\cup\{\infty\}$ be maximal such that
  $x_1=\dots=x_n=1$. If $n=\infty$, then
  \[\Theta(x)^f=\Theta(11\dots)^f=0^f=-1=\Theta(00\dots)=\Theta(x^{\tilde f});\]
  otherwise, $x=1\dots10x_{n+2}x_{n+3}\dots$; write
  $y=x_{n+3}x_{n+4}\dots$. If $x_{n+2}=0$, then
  \[\Theta(x)^f=(2^n+2^{n+1}+2^{n+2}\Theta(y))^f=2^{n+1}-1+2^{n+2}\Theta(y)=\Theta(0\dots001y)=\Theta(x^{\tilde f}),\]
  while if $x_1=1$, then
  \[\Theta(x)^f=(2^n+2^{n+2}\Theta(y))^f=2^{n+2}-1+2^{n+2}\Theta(y)=\Theta(0\dots000y)=\Theta(x^{\tilde f}).\]
\end{proof}

\subsection{Fibonacci sequence} From now on, we identify the semigroup
$F$ with $\Gamma(I)$ and use the notation $F$ for both.  Recall the
definition of the following elements of $F$: $f_1=s$, $f_2=f$, and
$f_n=f_{n-2}f_{n-1}$ for $n\ge3$, and the Fibonacci numbers $\Phi_n$
defined by $\Phi_1=\Phi_2=1$ and $\Phi_n=\Phi_{n-2}+\Phi_{n-1}$ for
$n\ge 3$. We note that $\|f_n\|=\Phi_n$ for all $n\ge1$. Set
$\varphi=(1+\sqrt5)/2$; recall that $\Phi^n\approx\varphi^n/\sqrt5$.

Note that if $n$ is even and $\ge6$, then $f_n$ starts with $fsfs$,
while if $n$ is odd and $\ge5$ then $f_n$ starts with $sffs$. The
following statements are easily proven by induction:
\begin{lem}\label{lem:veryelemrel}
  (1) For all $k\le n-2$ with $k\equiv n\pmod 2$ we have $f_n\equiv
  f_k(f_{k+1}f_{k+3}\cdots f_{n-1})$; in particular, $f_k$ is a prefix
  of $f_n$.

  (2) For all $k>1$ we have $f_{k+1}^2\equiv f_{k-1}f_{k+2}$.
\end{lem}


\subsection{Some relations} By~\eqref{eq:product}, we have
\begin{equation}\label{eq:prodzeta}
  \pair<g_0,g_1>\zeta\pair<h_0,h_1>\zeta=\pair<g_0h_0,g_1h_0>\zeta
\end{equation}
for all $g_i,h_i\in F$, a calculation that we will use repeatedly; so
by direct computation
\begin{lem}\label{lem:elemrel}
  In $F$ the following relations hold: $s^2=e$; $f^3=f$.
\end{lem}
\begin{proof}
  If $x=x_1\dots x_n$, then $x^{s^2}=((1-x_1)x_2\dots x_n)^s=x$, so
  $s^2=e$. Then
  \begin{align*}
    \phi(f^2)&=\pair<s,f>\zeta\pair<s,f>\zeta=\pair<s^2,fs>\zeta,\\
    \phi(f^3)&=\pair<s^2,fs>\zeta\pair<s,f>\zeta = \pair<s^3,fs^2>\zeta =
    \pair<s,f>\zeta=f. \qedhere
  \end{align*}
\end{proof}

\subsection{Contraction} The proofs will ultimately all rely on
some form of induction on the length of representations of semigroup
elements as words over $\{s,f\}$.

\begin{lem}\label{lem:contracting}
  The semigroup $F$ is contracting.
\end{lem}
\begin{proof}
  Consider $g\in F$, and write it as a word $w$ of minimal length
  $w_1\dots w_n$; write also $\phi(g)=\pair<g_0,g_1>\pi_g$. Then, by
  Lemma~\ref{lem:elemrel}, there cannot be more than two $f$'s in a
  row in $w$, so every group of three letters
  $w_{3k+1}w_{3k+2}w_{3k+3}$ contains at least an $s$, which will
  contribute no letter to $g_0$ nor $g_1$. The other two letters
  contribute at most one letter each to each of $g_0$ and $g_1$. In
  total, $\|g_0\|\le\frac23(\|g\|+2)$, and similarly for $\|g_1\|$.
\end{proof}
Note, in fact, that the contraction ratio $\eta$ may be chosen as
$1/\varphi$, with a little more care. This becomes apparent in the
next result.

\begin{lem}\label{lem:decompf_n}
  The decomposition of $f_n$ satisfies
  \begin{align*}
    \phi(f_1) &= \pair<e,e>\sigma, & \phi(f_2) &= \pair<s,f>\zeta,\\
    \phi(f_3) &= \pair<f,s>\zeta,  & \phi(f_4) &= \pair<f_3,f^2>\zeta,\\
    \phi(f_{2n}) &= \pair<f_{2n-1}, f_4f_6\cdots f_{2n-2}>\zeta\text{ if }n\ge 3,\\
    \phi(f_{2n+1}) &= \pair<f_{2n}, f_2(f_5f_7\cdots
    f_{2n-1})>\zeta\text{ if }n\ge 2.
  \end{align*}
\end{lem}
\begin{proof}
  The first four claims may be checked directly; we note that if
  $\phi(g)=\pair<g_0,g_1>\zeta$, then $\phi(sg)=\pair<g_1,g_0>\zeta$.
  Next, by induction,
  \begin{align*}
    \phi(f_{2n})&=\phi(f_{2n-2})\phi(f_{2n-1})
    =\pair<f_{2n-3},f_4f_6\cdots f_{2n-4}>\zeta\pair<f_{2n-2},*>\zeta\\
    &=\pair<f_{2n-3}f_{2n-2},f_4f_6\cdots f_{2n-4}f_{2n-2}>\zeta\\
    \intertext{as claimed, where $*$ stands for an irrelevant element
      of $F$; and}
    \phi(f_{2n+1})&=\phi(f_{2n-1})\phi(f_{2n})
    =\pair<f_{2n-2},f_2(f_5\cdots f_{2n-3})>\zeta\pair<f_{2n-1},*>\zeta\\
    &=\pair<f_{2n-2}f_{2n-1},f_2(f_5\cdots f_{2n-3}f_{2n-1})>\zeta.\qedhere
  \end{align*}
\end{proof}

\section{Proofs}\label{sec:proofs}
We start by finding relations in $F$; we then show that they yield a
simple normal form for elements of $F$.

\subsection{Relations}
Recall that for $n\ge1$ we defined
in~(\ref{eq:def:r_n},\ref{eq:def:r'_n}) words $r_n$ and $r'_n$ over
the alphabet $\{s,f\}$ by $r_1=s^2$, $r'_1=e$, and
\begin{align*}
  r_n &= f_{n+1}f_n^2 \equiv f_{n-1}f_n^3\;(\equiv
  f_{n+1}f_{n-2}f_{n+1}\text{ if }n\ge3),\\
  r'_n &= f_{n\%2+1}(f_{n\%2+5}f_{n\%2+7}\cdots f_{n-1}f_{n+1})f_n.
\end{align*}

Note that $r'_n$, for odd $n\ge 3$, can be obtained from $r_n$ by removing
the first four letters; indeed
\begin{equation}\label{eq:r_nodd}
  \begin{aligned}
    fsfsr'_n &\equiv f_4f_1\cdot f_2(f_6\cdots f_{n+1})f_n\equiv
    f_4f_3(f_6f_8\cdots f_{n+1})f_n\equiv f_4f_5^2(f_8\cdots f_{n+1})f_n\\
    &\equiv f_6f_5(f_8\cdots f_{n+1})f_n\equiv\cdots\equiv
    f_{n-3}f_{n-2}^2f_{n+1}f_n\equiv f_{n-1}f_n^3\equiv r_n.
  \end{aligned}
\end{equation}

Similarly, $r'_n$ for even $n\ge 4$ can be obtained from $r_n$ by
removing the first three letters and replacing them by $s$; indeed,
writing $r'_n=sr''_n$,
\begin{equation}\label{eq:r_neven}
  \begin{aligned}
    sffr''_n &\equiv f_3f_2(f_5\cdots f_{n+1})f_n\equiv
    f_3f_2(f_5f_7\cdots f_{n+1})f_n\equiv f_3f_4^2(f_7\cdots f_{n+1})f_n\\
    &\equiv f_5f_6^2(f_9\cdots f_{n+1})f_n\equiv\cdots\equiv
    f_{n-3}f_{n-2}^2f_{n+1}f_n\equiv f_{n-1}f_n^3\equiv r_n.
  \end{aligned}
\end{equation}

Let $e_n$ be the word obtained from $f_n$ by deleting its first two
symbols.
\begin{lem}\label{lem:delete2}
  If $n\ge 4$, then we have $\phi(f_n)=\pair<f_{n-1},e_{n-1}>\zeta$.
\end{lem}
\begin{proof}
  Assume first that $n$ is even. By Lemma~\ref{lem:veryelemrel} we
  have $f_{n-1}=f_3f_4f_6\dots f_{n-2}$, so $e_{n-1}=f_4f_6\cdots
  f_{n-2}$, and the Lemma holds by Lemma~\ref{lem:decompf_n}.
  Similarly, if $n$ is odd then $e_{n-1}=f_2f_5\cdots f_{n-2}$ by
  Lemma~\ref{lem:veryelemrel}, and again the Lemma holds by
  Lemma~\ref{lem:decompf_n}.
\end{proof}

\begin{lem}\label{lem:rel}
  In $F$ the relations $r_n=r'_n$ hold for all $n\ge1$.
\end{lem}
\begin{proof}
  The cases $n\le2$ are covered in~\ref{lem:elemrel}, since we have
  $r_2\equiv sf^3=sf\equiv r'_2$. Let us therefore assume $n\ge3$.  We
  follow the notation of Lemma~\ref{lem:delete2}.  We have the
  decomposition
  \[\phi(r_n) = \pair<f_nf_{n-1}^2,e_nf_{n-1}^2>\zeta =
  \pair<r_{n-1},fsr'_{n-1}>\zeta.\]
  Suppose first that $n$ is even.  Then by~\eqref{eq:prodzeta} we have
  \begin{align*}
    \phi(r'_n) &= \phi(f_1(f_5f_7\cdots f_{n+1})f_n) =
    \pair<e_4,f_4>\zeta\pair<f_6,*>\zeta\cdots\pair<f_n,*>
    \zeta\pair<f_{n-1},*>\zeta\\ 
    &= \pair<f_2(f_6f_8\cdots f_n)f_{n-1},f_4f_6\cdots f_nf_{n-1}>\zeta =
    \pair<r'_{n-1},fsr'_{n-1}>\zeta,
    \intertext{where $*$ stands for an element of $F$ that is not
      relevant to the calculation. If $n$ is odd, then similarly}
    \phi(r'_n) &= \phi(f_2f_6\cdots f_{n+1}f_n) =
    \pair<f_1,f_2>\zeta\pair<f_5,*>\zeta\cdots\pair<f_n,*>\zeta
    \pair<f_{n-1},*>\zeta\\ 
    &= \pair<f_1f_5\cdots f_nf_{n-1},f_2f_5\cdots f_nf_{n-1}>\zeta =
    \pair<r'_{n-1},fsr'_{n-1}>\zeta.
  \end{align*}
  We have $r_{n-1}=r'_{n-1}$ by induction, so $r_n=r'_n$.
\end{proof}

\subsubsection{More relations} We will ultimately show that
$\{r_n=r'_n\}$ is a complete set of relations for $F$; however, we
will first describe more relations, that are consequences of these but
allow much faster simplifications.

\begin{lem}\label{lem:abb}
  Consider $a,b\in\N$ with $a\ge3$ and $a\ge b+2$. Then
  \[f_af_b^2 = f_a.\]
\end{lem}
\begin{proof}
  We proceed by induction on $b$.  Since $a\ge b+2$, the word $f_a$
  ends with $f_{b+2}$; it therefore suffices to show that
  $f_{b+2}f_b^2 = f_{b+2}$.  We also note that the statement holds for
  $b\le2$, since then $f_3f_1^2 = f_3$ and $f_4f_2^2 = f_4$ by
  Lemma~\ref{lem:elemrel}. Then for $b\ge3$ we have
  \begin{align*}
    f_{b+2}f_b^2 &= f_bf_{b+1}f_b^2 = f_b r_{b} = f_b r'_{b} \\
    &= f_bf_{1+b\%2}f_{5+b\%2}f_{7+b\%2}\cdots f_{b-1}f_{b+1}f_b\\
    &= f_bf_{1+b\%2}^2f_{2+b\%2}f_{4+b\%2}f_{7+b\%2}f_{9+b\%2}\cdots f_{b-1}f_{b+1}f_b\\
    &= f_bf_{2+b\%2}^2f_{3+b\%2}f_{7+b\%2}f_{9+b\%2}\cdots f_{b-1}f_{b+1}f_b\\
    &= f_bf_{3+b\%2}f_{7+b\%2}f_{9+b\%2}\cdots f_{b-1} f_{b+1}f_b\\
    = \cdots &= f_bf_{b-3}f_{b+1}f_b = f_bf_{b-3}^2f_{b-2}^2f_{b-1}f_b
    = f_bf_{b-1}f_b = f_{b+2}.\qedhere
  \end{align*}
\end{proof}

\begin{lem}\label{lem:r'_n}
  For all $n\ge2$ we have
  \[r'_n = (f_{1}^2 f_{2}^2 f_{3}^2\cdots f_{n-3}^2 f_{n-2}^2)f_{n+1}.\]
\end{lem}
\begin{proof}
  If $n=2$ then $r'_2=f_3$ and the Lemma holds. If $n\ge3$, then
  \begin{align*}
    r'_n & = f_{1+n\%2} f_{5+n\%2} \cdots f_{n-1} f_{n+1} f_n\\
    & = f_{1+n\%2}^2 f_{2+n\%2}^2 f_{3+n\%2} f_{7+n\%2} \cdots
    f_{n-1} f_{n+1} f_n\\
    = \cdots & = f_{1+n\%2}^2 f_{2+n\%2}^2 f_{3+n\%2}^2 \cdots f_{n-3}^2 f_{n-2}^2 f_{n-1} f_n.\qedhere
  \end{align*}
\end{proof}

\begin{lem}\label{lem:fn 5=3}
  The equality $f_n^4 = f_n^2$ holds in $F$ for $n\le4$, and $f_n^5 =
  f_n^3$ holds in $F$ for all $n \ge 1$.
\end{lem}
\begin{proof}
  The first statement follows from the definition of $r_n$.  It
  suffices to check the second one for $n \ge 5$.  Using
  Lemmata~\ref{lem:r'_n} and~\ref{lem:abb}, we have
  \begin{align*}
    f_n^5 &\equiv f_n(f_{n-2} f_{n-1})f_{n}^3\equiv f_nf_{n-2}r_n =
    f_nf_{n-2} r'_{n} = f_n f_{n-2} f_{1}^2 f_{2}^2 f_{3}^2 \cdots
    f_{n-3}^2 f_{n-2}^2f_{n+1}\\
    &= \cdots = f_n f_{n-2} f_{n-3}^2 f_{n-2}^2 f_{n+1} =
    f_nr'_{n-3}f_{n-2}^2f_{n+1}
    = f_n f_{1}^2f_{2}^2 \cdots f_{n-5}^2 f_{n-2}^3 f_{n+1}\\
    &= \cdots = f_nf_{n-2}f_{n+1} = f_n f_{n-2} f_{n-1}f_n = f_n^3.\qedhere
  \end{align*}
\end{proof}


\begin{lem}\label{lem:reducer}
  For $n \ge 8$ we have
  \begin{align*} 
    f_{n-2} f_{n-3} f_{n-5} f_{n-7} \cdot f_{n-2} f_n &= f_{n-2}
    f_{n-3} f_{n-5} \cdot f_{n-3} f_n \\
    &= f_{n-2} f_{n-3} \cdot f_{n-4} f_n = f_{n-2} \cdot f_{n-5} f_n =
    f_n.
  \end{align*}
\end{lem}

\begin{proof}
  At first we prove by induction that for $n\ge1$ we have
  \[f_{n+1} f_n^2 f_{n+2} = f_{n+3}.\]
  The induction starts with two cases: for $n = 1$ we have $f_2 f_1^2
  f_3 = f_2 f_3 = f_4$, and for $n = 2$ we have $f_3 f_2^2 f_4 = f_3
  f_4 = f_5$.  Using twice Lemma~\ref{lem:r'_n} and the induction
  hypothesis, we have for $n \ge 3$
  \begin{align*}
    f_{n+1} f_n^2 f_{n+2} & = f_{n-1} f_{n}^3 f_n f_{n+1}\\
    & = f_{n-1} f_{n-2} \underbrace{f_{n-1} f_n^3}_{r'_n} f_{n+1}
    = f_{n-1} f_{n-2} f_1^2 f_2^2 \cdots f_{n-2}^2 f_{n+1}^2\\
    & = f_{n-1} \underbrace{f_{n-2} f_{n-3}^2}_{r'_{n-3}} f_{n-2}^2
    f_{n+1}^2 = f_{n-1} f_1^2 f_2^2 \cdots f_{n-5}^2 f_{n-2}^3
    f_{n+1}^2\\
    & = f_{n-1} f_{n-2}^3 f_{n+1}^2 = f_{n-1} f_{n-2}^2 \cdot
    f_{n-2} f_{n-1} f_{n} f_{n+1}\\
    & = f_{n-1} f_{n-2}^2 f_{n} \cdot f_{n+2} = f_{n+1} f_{n+2}
    = f_{n+3}.
  \end{align*}
  Using the arguments above, we carry out the following
  transformations
  \begin{align*}
    f_{n-2} & f_{n-3} f_{n-5} f_{n-7} \cdot f_{n-2} f_n =\\
    & = f_{n-2} \underbrace{f_{n-3} f_{n-5} f_{n-7} \cdot
      f_{n-6}}_{f_{n-3} f_{n-5}^2 = f_{n-3} } f_{n-5} f_{n-3} f_n\\
    & = f_{n-2} f_{n-3} f_{n-5} \cdot f_{n-3} f_m = f_{n-2} f_{n-3}
    f_{n-5} \cdot f_{n-5} f_{n-4} f_n \\
    & = f_{n-2} f_{n-3} \cdot f_{n-4} f_n = f_{n-2} f_{n-3} f_{n-4}^2 f_{n-3} f_{n-1}\\
    & = f_{n-2} \cdot f_{n-5} f_{n-4} f_{n-3} f_{n-1} = f_{n-2} \cdot f_{n-5} f_n\\
    & = f_{n-2} f_{n-3}^2 f_{n-1} = f_n.  \qedhere
  \end{align*}
\end{proof}

\subsection{A normal form for $F$}
We first show that every semigroup element can be written in the form
of Theorem~\ref{thm:nf}. We write $\mathcal N$ for the set of words
described there.

The words in $\mathcal N$ are those words that can be written as a
product $f_{i_1}\dots f_{i_m}\dots f_{i_n}$ where the sequence
$i_1,\dots,i_m$ is increasing by steps of at least two, and the
sequence $i_m,\dots,i_n$ is decreasing by steps of at least one.

An arbitrary word $w\in\{s,f\}^*$ is put into normal form by replacing
subwords that are left-hand sides of the rules (N1)-(N10) by their
right-hand side, until no such substitution is possible. For example,
consider the possible reductions that can be performed on
the word $f_4f_3^2f_4$:
\[\xymatrix{
  {f_4f_3^2f_4}\ar[rr]^{\text{N2}:f_3f_4\to f_5}\ar[d]_{\text{N3}:f_3^2\to f_1f_4}
  && {f_4f_3f_5}\ar[rr]^{\text{N8}:p=1,a=3,q=2}
  && {f_2f_5}\ar[rr]^{\text{N5}:f_2\to f_1f_3} && {f_1f_3f_5.}\\
  {f_4f_1f_4^2}\ar[rr]^{\text{N3}:f_4^2\to f_2f_5} \ar@/_1pc/[rrrr]_{\text{N10}:p=3,a=1,q=3}
    && {f_4f_1f_2f_5}\ar[u]_{\text{N2}:f_1f_2\to f_3}
    && {f_4^2}\ar[u]_{\text{N3}:f_4^2\to f_2f_5}
  }
\]

\begin{lem}\label{lem:terminating}
  The rewriting system in Theorem~\ref{thm:rws} is terminating.
\end{lem}
\begin{proof}
  We define a function $\eta=(\eta_1,\eta_2):\{f_1,f_2,\dots\}^*\to\N^2$
  with the lexicographic ordering, such that if $w\to w'$ is obtained
  by an elementary reduction then $\eta(w)>\eta(w')$. Since $\N^2$ is
  well ordered this will prove that the rewriting system is
  terminating.

  Consider therefore $w\in\{f_1,f_2,\dots\}^*$, say $w=f_{i_1}\dots
  f_{i_n}$. Define
  \begin{align*}
    \eta_1(w)&=\begin{cases}
      \sum_{j=1}^n\Phi_{i_j}-3 & \text{ if }n\ge2\text{ and }i_1=1\text{
        and }i_2\ge3\text{ is odd},\\
      \sum_{j=1}^n\Phi_{i_j} & \text{ otherwise};
    \end{cases}\\
    \eta_2(w)&=\sum_{j=1}^n(n+1-j)i_j.
  \end{align*}
  We may then check that $\eta_1(w)>\eta_1(w')$ if $w'$ is obtained
  from $w$ using one of the rules (N1,4--5,8--10), while
  $\eta_1(w)=\eta_1(w')$ and $\eta_2(w)>\eta_2(w')$ if $w'$ is
  obtained from $w$ using one of the rules (N2--3,6--7).

  For example, in rule (N9), only relation $r_{a+2}\to r'_{a+2}$ is
  used, and it shortens the length (as a word over the alphabet
  $\{s,f\}$) of $w$ by $4$; note that the length of $w$ is precisely
  $\sum\Phi_{i_j}$.  Write $w'=f_{i'_1}\dots f_{i'_{n'}}$; then
  \[\eta_1(w)\ge\sum_{j=1}^n\Phi_{i_j}-3>\sum_{j=1}^n\Phi_{i_j}-4=\sum_{j=1}^{n'}\Phi_{i'_j}\ge\eta_1(w').\]

  Consider, as another example, rule (N6). We have
  $\Phi_{a+1}+\Phi_a+\Phi_{a+3}=\Phi_{a+3}+\Phi_{a+2}$, so
  $\eta_1(w)=\eta_1(w')=\sum\Phi_{i_j}$. Say that the left-hand side
  of the rule appears at positions $n-j$, $n-j+1$ and $n-j+2$ in $w$;
  then
  \begin{align*}
    \eta_2(w)-\eta_2(w')&\ge(j+1)(a+1)+ja+(j-1)(a+3)-j(a+3)-(j-1)(a+2)\\
    &=(j+1)(a-1)+1>0.
  \end{align*}
  All other rules are proven to decrease $\eta$ in a similar fashion.
\end{proof}

The following calculations will be useful in the proof of
Lemma~\ref{lem:reddiff}. Let $u = f_{i_1} \dots f_{i_m} \dots f_{i_n}$
be an arbitrary semigroup word written in normal form such that $n \ge
1$, and $i_1, i_n > 1$, and $i_m$ is the maximal index. If $i_1=1$ or
$i_n=1$, the argument below should be applied to the word obtained
from $u$ by removing these instances of $f_1$, adapting the statements
accordingly. If $i_1=3$, let $k$ be the largest integer such that $i_r
= 2r+1$ for all $r < k$; if $i_1=4$, let $k$ be the largest integer
such that $i_r = 2r+2$ for all $r<k$; if $i_1\ge5$, then $k$ need not
be defined.  Then, by Lemma~\ref{lem:decompf_n}, the element $u$
affords the decomposition $\phi(u) = \pair<u_0,u_1>\zeta$, where
\begin{align}
 \label{eq:u0}
  u_0 &= \begin{cases}
    f_{i_1-1} \cdots f_{i_m-1} \cdots f_{i_n-1} & \text{ if $i_1 > 3$,}\\
    f_1 f_{2k-1} f_{i_k-1} \cdots f_{i_n-1} & \text{ if $i_1 = 3$,} \\
  \end{cases}\\
 \label{eq:u1}
  u_1 &= \begin{cases}
    (f_1 f_3 \cdots f_{i_1-2}) f_{i_2-1} \cdots f_{i_n-1}
    & \text{ if $i_1$ is odd and $i_2 \neq i_1-1$,} \\
    (f_4 f_6 \cdots f_{i_1-1}) f_{i_3-1} \cdots f_{i_n-1}
    & \text{ if $i_1$ is odd and $i_2 = i_1-1$,} \\
    (f_4 f_6 \cdots f_{i_1-2}) f_{i_2-1} \cdots f_{i_n-1}
    & \text{ if $i_1$ is even $\ge 6$ and $i_2 \neq i_1-1$,} \\
    (f_1 f_3 \cdots f_{i_1-1}) f_{i_3-1} \cdots f_{i_n-1}
    & \text{ if $i_1$ is even $\ge 6$ and $i_2 = i_1-1$,} \\
    f_1 f_{2k-1} f_{2k-2} f_{i_k-1} \cdots f_{i_n-1}
    & \text{ if $i_1 = 4$, $i_k \le 2k-2$,} \\
    f_{i_2-1} f_{i_3-1} \cdots f_{i_n-1}
    & \text{ if $i_1 = 4$, $i_k = 2k-1$ or $i_k \ge 2k+3$.}
  \end{cases}
\end{align}
The maximal index of $u_1$ is always less than $i_m$, but the maximal
index of $u_0$ is equal to $i_m$ if and only if $i_1 = 3$ and $i_k
\le 2k-2$.

As an example, consider $u=f_1f_3f_8$ and $v=f_3f_5f_7$. For $u$, we
first compute $\phi(f_3f_8)=\pair<f_1f_3f_7,f_1f_7>\zeta$ using the
first lines in~\eqref{eq:u0} and~\eqref{eq:u1} and $k=2$; so
\[\phi(u)=\pair<f_1f_7,f_1f_3f_7>\zeta.\qquad\text{Similarly,
}\phi(v)=\pair<f_1f_7,f_1f_4f_6>\zeta,\]
using the first lines in~\eqref{eq:u0} and~\eqref{eq:u1} and $k=4$. Note
that $u_0=v_0$ while $u\neq v$; this happened because one of $u,v$
started with $f_1$.

\begin{lem}\label{lem:2eq}
  The equations $sfsf^2x=fx$ and $sfy=fsy$ have no solution in $F$.
\end{lem}
\begin{proof}
  We will show that if the first equation had a solution $x$, then the
  second equation would have a solution $y$ that is not longer than
  $x$; and if the second equation had a solution $y$, then the first
  equation would have a solution $x$ that is shorter than $y$.

  Assume therefore that $x$ is a solution to the first equation, and
  write $\phi(x)=\pair<x_0,x_1>\pi$. Then
  $\phi(sfsf^2x)=\pair<f^2sx_0,sfsx_0>\zeta\pi=\pair<sx_0,fx_0>\zeta\pi$,
  so in particular, by considering the second co\"ordinate, $y=sx_0$
  is a solution to the second equation, and is not longer than $x$.

  Assume then that $y$ is a solution to the second equation, and write
  $\phi(y)=\pair<y_0,y_1>\rho$.  Then
  $\pair<fy_0,sy_0>\zeta\rho=\pair<sy_1,fy_1>\zeta\sigma\rho$, so
  $y_0=sfy_1=sfs(sy_1)=sfsfy_0$ and similarly $y_1=sfsfy_1$. Clearly
  $y\not\in\langle s\rangle$; so $y=s^ifz$ for some $i\in\N$ and $z\in
  F$. If $i$ is odd then $y_0$ begins with $f$, while if $i$ is even
  then $y_1$ begins with $f$. In all cases one of $y_0,y_1$ begins
  with $f$, say $y_i=fx$; then $x$ is a solution to the first
  equation, and is shorter than $y$.
\end{proof}

\begin{lem}\label{lem:reddiff}
  All elements of $\mathcal{N}$ are distinct in $F$.
\end{lem}
\begin{proof}
  The elements $e$ and $f_1$ are distinct since they act differently
  on $X$. In addition, they differ from all other elements of $F$ since
  they are invertible.

  Consider now two distinct normal forms
  \begin{xalignat*}{2}
    u &= f_1^\delta f_{i_1}\cdots f_{i_m}\cdots f_{i_n},&
    v &= f_1^\epsilon f_{j_1}\cdots f_{j_p}\cdots f_{j_q},
  \end{xalignat*}
  with $i_m, j_p \ge 3$ and $i_m \ge j_p$.

  If $i_n = 1$ and $j_q \neq 1$, then $0^u = 1$ while $0^v = 0$ so
  they are distinct; if $i_n \neq 1$ and $j_q = 1$ the same argument
  applies; if $i_n = j_q = 1$ then we may cancel both from the normal
  form and proceed. Let us therefore suppose $i_n, j_q > 1$.

  We proceed by induction on the maximal index $i_m$, to prove that
  all normal forms with at most that maximal index are distinct in
  $F$.

  The induction starts with the elements $f_3$, $f_3 f_2$, $f_1 f_3$
  and $f_1 f_3 f_2$. They afford the following decompositions:
  \begin{align*}
    \phi(f_3) & = \pair< f_1 f_3, f_1 >\zeta, &
    \phi(f_3f_2) & = \pair< f_1 f_3 f_1, e >\zeta,\\
    \phi(f_1 f_3) & = \pair< f_1, f_1 f_3 >\zeta, &
    \phi(f_1f_3f_2) & = \pair< e, f_1 f_3 f_1 >\zeta.
  \end{align*}
  The first co\"ordinates of these decompositions are different by the
  remarks above.  Let us therefore suppose $i_m \ge 4$.

  \begin{description}
  \item[If $\boldsymbol{\delta=\epsilon=1}$] up to left-multiplying
    $u$ and $v$ by $f_1$, we may replace $\delta$ and $\epsilon$ by
    $0$ and proceed.
  \item[If $\boldsymbol{\delta=\epsilon=0}$] then
    \begin{align*}
      \phi(u) &= \pair<f_{i_1-1} \cdots f_{i_m-1} \cdots f_{i_n-1}, u_1 >\zeta,\\
      \phi(v) &= \pair<f_{j_1-1} \cdots f_{j_p-1} \cdots f_{j_q-1},
      v_1 >\zeta,
    \end{align*}
    where $u_1$ and $v_1$ are defined by~\eqref{eq:u1}.

    If $i_1, j_1 > 3$, then the first co\"ordinates of $\phi(u)$ and
    $\phi(v)$ are in normal form, and are distinct by assumption, so
    by induction they are distinct in $F$. Similarly, if $i_1 = j_1 =
    3$ then the second co\"ordinates equal $f_1 f_{i_2-1} \cdots
    f_{i_n-1}$ and $f_1 f_{j_2-1} \cdots f_{j_p-1} \cdots f_{j_q-1}$,
    and they are in normal form, so by induction $u$ and $v$ are
    distinct in $F$.

    We consider therefore the case $i_1 = 3,j_1>3$. Let $k$ be maximal
    such that $i_r=2r+1$ for all $r<k$. Then the first co\"ordinate of
    $\phi(v)$ is in normal form and starts by $f_{j_1-1}$ while the
    first co\"ordinate of $\phi(u)$, when put in normal form, is $f_1
    f_{2k-1} f_{i_k-1} \cdots f_{i_n-1}$. By the next item or by the
    induction hypothesis these elements are distinct in $F$, whence
    $u$ and $v$ are distinct in $F$, too.

    The case $i_1 = 3, j_1 > 3$ is considered similarly.

  \item[If $\boldsymbol{\delta\neq\epsilon}$] by symmetry, we may
    assume $\delta=0$ and $\epsilon=1$. Then
    \[
    \phi(v) = \pair< e_{j_1-1} f_{j_2-1}\cdots f_{j_q-1}, f_{j_1-1}
    \cdots f_{j_p-1} \cdots f_{j_q-1} > \zeta.
    \]
    Assume for contradiction that $u,v$ are equal in $F$. Then, in
    $F$, we have
    \begin{align*}
        u_0& =f_{i_1-1}\cdots f_{i_n-1} & &= & v_1 &= e_{j_1-1} f_{j_2-1}\cdots
        f_{j_q-1},\\
        u_1 &= e_{i_1-1} f_{i_2-1}\cdots f_{i_n-1} & &= & v_0 &= f_{j_1-1}\cdots
        f_{j_q-1},
    \end{align*}
    so these words, when put in normal form, must be equal. It follows
    from Lemma~\ref{lem:delete2} that $u_0 = sf u_1$ if $i_1$ is odd
    and $u_0 = fs u_1$ if $i_1$ is even. Similarly, we have $v_0 = sf
    v_1$ or $v_0 = fs v_1$, depending on the parity of $j_1$.

    We consider in turn all possibilities for the parity of $i_1$ and
    $j_1$.  If $i_1$ is even and $j_1$ is odd, then $u_0=sf v_0=sffs
    v_1 = v_1$; it follows from~\eqref{eq:u1} that when put in normal
    form, $v_1$ starts with $f$.  Using Lemma~\ref{lem:2eq}, we reach
    the contradiction. The same argument holds, by symmetry, if $i_1$
    is odd and $j_1$ is even.


    If $i_1$ and $j_1$ are both even, then we find $(sf)^2u_0=u_0$; now
    by assumption $u_0$ is of the form $f_k\dots$ for an odd $k$, so
    its normal form starts with $f$; this contradicts Lemma~\ref{lem:2eq}.

    Finally, if $i_1$ and $j_1$ are both odd, then we find
    $(fs)^2u_0=u_0$, and $u_0$ is of the form $f_k\dots$ for an even
    $k$, so its normal form starts with $s$; again this contradicts
    Lemma~\ref{lem:2eq}.
  \end{description}
\end{proof}

\begin{proof}[Proof of Theorem~\ref{thm:rws}]
  We saw in Lemma~\ref{lem:terminating} that the rewriting system is
  terminating. On the other hand, consider a word
  $w\in\{f_1,f_2,\dots\}^*$ and two reduced words $w',w''$ that are
  gotten from $w$ by applying two different sequences of elementary
  reductions. The reduced forms of the rewriting system correspond
  precisely to the normal form described in Theorem~\ref{thm:nf}; and
  they define equal elements in $F$, so by Lemma~\ref{lem:reddiff}
  they are graphically equal. This proves that the rewriting system is
  confluent.
\end{proof}

\begin{proof}[Proof of Theorem~\ref{thm:nf}]
  Every semigroup element can be represented by a word, and this word
  can be put in reduced form by applying the rules from
  Theorem~\ref{thm:rws}. As noted above, this reduced form coincides
  with the normal form; therefore the natural map $\mathcal{N}\to F$
  is onto. It is one-to-one by Lemma~\ref{lem:reddiff}.

  Finally, all the rules (N1-10) derive from the rules $s^2=1$ (N1)
  and $r_n\to r'_n$ for $n\ge2$ (N4,8-10), and from graphical
  equalities (N2-3,5-7). Each of the rules $r_n\mapsto r'_n$ is
  strictly length-shortening, so $\mathcal{N}$ is precisely the set of
  irreducible words with respect to the rules, and consists of minimal
  representative words. In other words, if there were a shorter
  representation $w'$ of some $w\in\mathcal N$, then this shorter
  representation could be reduced to a word $w''\in\mathcal N$ using
  the rules $r_n\mapsto r'_n$, still yielding a word shorter, and
  therefore different, from $w$. This contradicts
  Lemma~\ref{lem:reddiff}.

  Now $\|f_n\|=\Phi_n$, and the initial $f_1$ should cancel with the
  initial $s$ of $f_{i_1}$ if $i_1$ is odd; this explains the length
  formula.
\end{proof}

\begin{proof}[Proof of Theorem~\ref{thm:pres}]
  The relations $r'_n=r_n$ for $n\ge1$ hold in $F$ by
  Lemma~\ref{lem:rel}. They are also sufficient; otherwise, there
  would be two distinct words $w,w'$ that are unequal in the abstract
  semigroup $\langle s,f|\,r'_n=r_n\rangle$ but are equal in
  $F$; by Theorem~\ref{thm:nf} the normal form representation of $w$
  and $w'$ are equal; but then $w$ and $w'$ can be transformed into
  each other using only the relations $r_n=r'_n$, contradicting our
  initial assumption.
\end{proof}

\subsection{Growth}
\begin{proof}[Proof of Theorem~\ref{thm:growth}]
  Consider a typical element in the ball of radius $\ell$ in $F$. It will
  have a minimal normal form as given by Theorem~\ref{thm:nf}, whence
  $\Phi_{i_m}\le\ell$. The maximal length of a normal form for that
  value of $i_m$ is
  $\Phi_{i_m}+\Phi_{i_m-1}+2\Phi_{i_m-2}+\Phi_{i_m-3}+2\Phi_{i_m-4}+\cdots
  < 2\Phi_{i_m+1}$, so it follows that $\ell<2\Phi_{i_m+1}$.  Now since
  $\Phi_i\approx\varphi^i/\sqrt 5$, we obtain
  \[\varphi^{i_m}/\sqrt 5 \le \ell < 2/\sqrt 5\varphi^{i_m+1},\]
  so
  \[\frac{\log(\sqrt5\ell/(2\phi))}{\log\phi}\le i_m\le\frac{\log(\sqrt5\ell)}{\log\phi}.\]
  Given any value of $i_m$, there are $2^{i_m}$ choices of $i_j$ to
  include in the normal form with $j>m$, and $2\Phi_{i_m-2}$ choices
  for the $i_j$ to include with $j<m$: the $2$ counts the possible
  values of $i_1$, and the $\Phi_{i_m-2}$ counts the possible values
  of $(i_3,\dots,i_{m-2})\in\{0,1\}^{m-4}$ with no two consecutive
  $1$'s. It follows that the ball growth function $\gamma(\ell)$ is
  given by
  \[\gamma(\ell)=\sum_{k=1}^{i_m}2^k2\Phi_{k-2}\approx \frac2{\sqrt
    5\varphi^2}\frac{(2\varphi)^{i_m}-1}{2\varphi-1},
  \]
  which gives the upper and lower bounds
  \[\frac{2}{\sqrt
  5\varphi^2(2\varphi-1)}\left(\frac{\sqrt5\ell}{2\varphi}\right)^\alpha\le
  \gamma(\ell)\le\frac{2}{\sqrt 5\varphi^2(2\varphi-1)}\big(\sqrt5\ell\big)^\alpha,
  \]
  with $\alpha=\log(2\varphi)/\log\varphi$.
\end{proof}

\subsection{Identities}
\begin{proof}[Proof of Theorem~\ref{thm:identity}]
  More generally, we will show that all of the equalities
  $g(wg)^5=g(wg)^3$, for all $w\in\{1,s,f,fs,f^2,f^2s\}$, hold in $F$.
  We prove this by considering simultaneously all these equalities,
  and proceed by induction on $\|g\|$.

  The induction is easily checked if $\|g\|\le2$, by considering in
  turn the cases $g=1,s,f,sf,fs,f^2$.  For longer $g$, depending on
  whether $g$ starts or ends with an $s$, we have
  $\phi(g)=\pair<fh,sh>\zeta$ or $\pair<sh,fh>\zeta$ or
  $\pair<fh,sh>\zeta\sigma$ or $\pair<sh,fh>\zeta\sigma$; for
  simplicity, assume that we are in the first case. Then $\|h\|<\|g\|$
  by Lemma~\ref{lem:contracting}, so the equalities above hold for $h$
  by induction. Then
  \begin{align*}
    \phi(gg^5) &= \pair<fh(fh)^5,sh(fh)^5>\zeta = \phi(gg^3),\\
    \phi(g(sg)^5) &= \pair<fh(sh)^5,sh(sh)^5>\zeta = \phi(g(sg)^3),\\
    \phi(g(fg)^5) &= \pair<fh(f^2h)^5,sh(f^2h)^5>\zeta = \phi(g(fg)^3),\\
    \phi(g(fsg)^5) &= \pair<fh(fsh)^5,sh(fsh)^5>\zeta = \phi(g(fsg)^3),\\
    \phi(g(f^2g)^5) &= \pair<fh(f^3h)^5,sh(f^3h)^5>\zeta =
    \pair<fh(fh)^5,sh(fh)^5>\zeta = \phi(g(f^2g)^3),\\
    \phi(g(f^2sg)^5) &= \pair<fh(f^2sh)^5,sh(f^2sh)^5>\zeta = \phi(g(f^2sg)^3).\qedhere
  \end{align*}
\end{proof}

\subsection{Ideal structure}
\begin{proof}[Proof of Theorem~\ref{thm:quotients}]
  Let us prove by induction on $n$ that $z_n(u) = 0^n$ for any $u \in
  X^n$, and that $f_{n+2} f_{n+1} = z_n$ holds in $W_n$, i.e.\ that
  $z_n$ and $f_{n+2} f_{n+1}$ are the left-hand zeroes over $X^n$. For
  $n = 1$ it follows from the definition of $f_3$ that $f_3 : X^1 \to
  \{0\}$, whence the relations $f_1 f_3 = f_2 f_3 = f_3$ and $f_3 f_2
  = f_3$ hold.

  Choose $n \ge 2$. The element $z_n$ has the following decomposition:
  \[
  z_n = \pair<f_2 f_4 \cdots f_{n+1}, f_1 f_4 f_6 \cdots f_{n+1}>
  \zeta = \pair< f_1 f_3 z_{n-1}, f_1 z_{n-1}> \zeta
  \]
  if $n$ is odd, and
  \[
  z_n = \pair<f_3 f_5 \cdots f_{n+1}, f_2 f_2 f_5 f_7 \cdots f_{n+1}> \zeta =
  \pair<z_{n-1}, f_1 f_3 f_1 z_{n-1}> \zeta
  \]
  if $n$ is even. Clearly both co\"ordinates of decompositions end with
  $z_{n-1}$, which is a left-hand zero over $X^{n-1}$ by the
  induction hypothesis. As $z_n(x) = 0$ for $x \in X$, we have $z_n :
  X^n\to\{0^n\}$.

  Similarly, the element $f_{n+2} f_{n+1}$ is written in the
  following way:
  \begin{multline*}
    f_{n+2} f_{n+1} = \pair<f_{n+1} f_n, f_1 f_3 \cdots f_{n} f_{n}>
    \zeta \\= \pair< f_{n+1} f_n, f_4 \cdots f_{n-1} f_{n+1}> \zeta =
    \pair< f_{n+1} f_n, z_{n-1}> \zeta
  \end{multline*}
  if $n$ is odd, and
  \begin{multline*}
    f_{n+2} f_{n+1} = \pair<f_{n+1} f_n, f_4 f_6 \cdots f_{n}
    f_{n}> \zeta \\= \pair< f_{n+1} f_n, f_1 f_3 \cdots f_{n-1}
    f_{n+1}> \zeta = \pair< f_{n+1} f_n, f_1 z_{n-1}> \zeta
  \end{multline*}
  if $n$ is even. By the induction hypothesis and the proof above for
  $z_n$, the element $f_{n+2} f_{n+1}$ is a left-hand zero over
  $X^n$.

  As $z_n, f_{n+2} f_{n+1} : X^n \to\{0^n\}$, the relations $f_1 z
  = f_2 z = z$ and $f_{n+2} f_{n+1} = z$ hold in $W_n$. We show by
  induction on the maximal index of a semigroup element that $W_n$
  includes no elements with maximal index $> n+2$. It follows from
  the definition of $f_n$ that
  \[
  f_{n+3} = f_{3-n\%2} f_{4-n\%2} f_{6-n\%2} \cdots f_{n+2}
  = f_{3-n\%2} z_n = z_n.
  \]
  For any $k > n+3$ the equality $f_k = f_{k-2} f_{k-1}$ holds
  and by the induction hypothesis $f_k = f_{k-2} z_n = z_n$.

  It is necessary to prove that the other semigroup elements define
  different transformations of $X^n$. Once again, we use induction on
  $n$. For $n = 1$ the four semigroup elements $e$, $f_1$, $f_3$ and
  $f_3 f_1$ define different automatic transformations over $X^1$.

  Fix $n \ge 2$. Let us write
  \begin{align*}
    S_n & = \{ g \in F \text{ with maximal index } n\},\\
    Z_n & = \{ g \in S_{n+2}:\; g \text{ includes }z_n\},\\
    A_n & = \{ g \in S_{n+2}:\; g \text{ includes neither $z_n$ nor
      $f_{n+2} f_{n+1}$} \},\\
    B_n & = S_{n+2} \setminus(Z_n\cup A_n).
  \end{align*}
  The semigroup $W_n$ is constructed as $W_{n-1} \cup B_{n-1} \cup Z_n
  \cup A_n$. It follows from the proof above that elements with
  maximal index $n+2$ act as left-hand zeroes on $X^{n-1}$. Then
  elements of $Z_n$ and $A_n$ coincide with elements $Z_{n-1}$ over
  $X^{n-1}$. On the other hand, let $g$ be an arbitrary element from
  $S_{n+2}$. It follows from~\eqref{eq:u1} that any input symbol
  decreases the maximal index of at least one co\"ordinate in the
  decomposition of $g$. Hence, there exists the input word $u$ of
  length $n-1$ such that it transforms $g$ into the element with
  maximal index $3$. As elements from $S_3$ can be distinguished over
  $X^2$, elements from $S_{n+2}$ cannot act as a left-hand zero over
  $X^{n+1}$. Therefore elements from $Z_n$ can be distinguished from
  elements of $B_{n-1} \cup W_{n-1}$ over $X^n$.

  Let $g$ be an arbitrary element from $A_n$. As $g$ does not include
  $z_n$ and $f_{n+2} f_{n+1}$, the second (or the first, if $g$ starts
  with $f_1$) co\"ordinate of its decomposition either contains
  $z_{n-2}$ and has maximal index $n$, or it does not include
  $z_{n-1}$ and has maximal index $n + 1$. In both cases the second
  co\"ordinate does not act as a left-hand zero over $X^{n-1}$ and
  therefore $g$ isn't a left-hand zero over $X^n$. Consequently, $A_n
  \cap Z_n = \emptyset$, and all elements of $W_n$ define different
  transformations of $X^n$.

  Now we count the elements in $W_n$. An arbitrary element with
  maximal index $k \ge 3$ can be written as
  \[
  f_1^{\alpha_1} \underbrace{f_3^{\alpha_3} \cdots
    f_{k-2}^{\alpha_{k-2}}}_{\Phi_{k-2} \text{ possibilities}} f_k
  \underbrace{f_{k-1}^{\beta_{k-1}} \cdots f_2^{\beta_2}
    f_{1}^{\beta_{1}}}_{2^{k-1} \text{ possibilities}},
  \]
  where $\beta_j \in\{ 0, 1\}$ for all $j = 1, 2, \dots, k-1$, and
  $\alpha_i \alpha_{i+1} = 0$, $i = 3, \dots, k-2$. Therefore there
  are $2\cdot\Phi_{k-2}\cdot 2^{k-1}$ elements with maximal index
  $k$. The semigroup $W_n$ includes all elements with maximal index
  that is $< n+2$, semigroup elements of the form $z
  f_{k-1}^{\beta_{k-1}} \dots f_2^{\beta_2} f_{1}^{\beta_{1}}$, and
  semigroup elements with maximal index $(n+2)$ that include neither
  $z_n$ nor $f_{n+2} f_{n+1}$.  Therefore,
  \begin{align*}
    \#W_n &= 2+\sum_{k = 3}^{n+1} 2^{k} \Phi_{k-2}+2^n+\left(
      \frac{1}{2} 2^{n+2} \Phi_{n}-2^{n+1} \right) \\
    &= \sum_{k = 1}^n 2^{k+2} \Phi_{k}-2^n-2^{n+1} \Phi_{n}+2.\qedhere
  \end{align*}
\end{proof}

To prove Corollary~\ref{cor:trace}, we begin by a Lemma:
\begin{lem}\label{lem:idzn} In the semigroup $F$, we have
  \begin{align}
    f_nf_{n-1} &= \begin{cases}
      fsfsf_{n+1} & \text{ if $n\ge 4$ is even},\\
      sffsf_{n+1} & \text{ if $n\ge 5$ is odd};
    \end{cases}\\
    z_n &= \begin{cases}
      fsf_{n+3} & \text{ if $n\ge 2$ is even},\\
      sfsf_{n+3} & \text{ if $n\ge 3$ is odd};
    \end{cases}\\
    f_n &= \begin{cases}
      fsfsf_{n-1}f_{n-2} = fz_{n-3} & \text{ if $n\ge 6$ is even},\\
      sffsf_{n-1}f_{n-2} = sfz_{n-3} & \text{ if $n\ge 5$ is odd};
    \end{cases}
  \end{align}
\end{lem}
\begin{proof}
  Induction on $n$, the basis of the induction being easily checked by
  applying the reduction algorithm. For example, take $n\ge 6$ even;
  then $fsfsf_{n+1}=fsfsf_{n-1}f_n=f_{n-2}f_{n-3}f_n$ by induction.
  Then
  $f_{n-2}f_{n-3}f_n=f_{n-2}f_{n-3}f_{n-2}f_{n-1}=f_{n-2}f_{n-1}^2=f_nf_{n-1}$.
  All other cases are treated similarly.
\end{proof}

\begin{proof}[Proof of Corollary~\ref{cor:trace}]
  The proof follows from Theorem~\ref{thm:quotients}. Assume first
  that $g$ includes $z_{n-2}$ or $f_n f_{n-1}$. Then it acts as a
  left-hand zero on $X^{n-2}$, whence $\tau(g) \le 2^{-(n-2)}$. On
  the other hand, $g$ belongs to $W_{n-1}$ and doesn't act as a
  left-hand zero on $X^{n+1}$; therefore $\tau(s) = 2^{-(n-2)}$.

  Now assume $g$ includes neither $z_{n-2}$ nor $f_n f_{n-1}$.
  Then $g$ belongs to $W_{n-2}$, but acts as a left-hand zero on
  $X^{n-3}$.  Therefore $g(X^*) = x_{\beta_1} x_{\beta_2} \dots
  x_{\beta_{n-3}} X^*$, whence $\tau(g) = 2^{-( n-3)}$.

  Note next that the generator of the ideal $I_n$ has trace $2^{-n}$.
  By the proof above, every element of trace $2^{-n}$ is a multiple
  $f_{n+3}$, of $f_{n+2}f_{n+1}$, or of $z_n$. By
  Lemma~\ref{lem:idzn}, the ideals generated by these three elements
  coincide.
\end{proof}

\begin{proof}[Proof of Theorem~\ref{thm:ideals}]
  Let $g$ be an arbitrary element with maximal index $n \ge 3$,
  generating the ideal $I = FgF$. We show that there exist $\ell(g),
  r(g) \in F$ such that
  \[
  \ell(g) \cdot g \cdot r(g) = \begin{cases} f_{n+1} & \text{ if $g$
      includes $z_{n-2}$ or $f_n f_{n-1}$ in
      its normal form,} \\
    f_n & \text{ otherwise.}
  \end{cases}
  \]
  Therefore $I = F g F = F f_k F = F_{2^{-k+3}}$ by
  Corollary~\ref{cor:trace}, where $k = n$ or $k = n+1$ depending on
  the cases above.

  Let us assume that $g = f_n f_{n-1}^{\beta_{n-1}} \cdots
  f_{2}^{\beta_2} f_1^{\beta_1}$, where $\beta_i \in \{ 0,1 \}$. Then
  we claim that for the element $g_2 = f_1^{\beta_1} f_{2}^{\beta_2}
  \cdots f_{n-2}^{\beta_{n-2}}$ the equality
  \begin{equation}\label{eq:reducer_right}
    g \cdot g_2 = f_n f_{n-1}^{\beta_{n-1}} \cdots f_{2}^{\beta_2} f_1^{\beta_1} \cdot
    f_1^{\beta_1} f_{2}^{\beta_2} \cdots f_{n-2}^{\beta_{n-2}}
    = f_n f_{n-1}^{\beta_{n-1}}
  \end{equation}
  holds. Applying $r_2=r_2'$ or Lemma~\ref{lem:abb} if
  necessary,
  \[
  f_n f_{n-1}^{\beta_{n-1}} \cdots f_{2}^{\beta_2} f_1^{\beta_1} \cdot
  f_1^{\beta_1} f_{2}^{\beta_2} = f_n f_{n-1}^{\beta_{n-1}} \cdots
  f_{3}^{\beta_3}.
  \]
  If $\beta_i = 0$, $\beta_{i-1} = 1$ for some $i \ge 4$ then it
  follows from Lemma~\ref{lem:abb} that
  \[
  f_n f_{n-1}^{\beta_{n-1}} \cdots f_{i+1}^{\beta_{i+1}} f_{i}^{0}
  f_{i-1} \cdot f_{i-1} f_i^0 = f_n f_{n-1}^{\beta_{n-1}} \cdots
  f_{i+1}^{\beta_{i+1}}.
  \]
  Otherwise, if $\beta_i = \beta_{i-1} = 1$ for some $i \ge 4$ then we
  have
  \[
  f_n f_{n-1}^{\beta_{n-1}} \cdots f_{i+1}^{\beta_{i+1}}
  f_{i} f_{i-1} \cdot f_{i-1} = f_n f_{n-1}^{\beta_{n-1}}
  \cdots f_{i+1}^{\beta_{i+1}} f_i
  \]
  by applying the relation
  \[
  f_p f_{m} f_{m-1} \cdot f_{m-1} = f_p f_1^2 f_2^2
  \cdots f_{m-3}^2 f_m = f_p f_m
  \]
  which holds for all $p>m\ge2$, due to Lemmata~\ref{lem:r'_n}
  and~\ref{lem:abb}.

  Now let $g = f_{3}^{\alpha_3} \cdots f_{n-2}^{\alpha_{n-2}}
  f_n$ with $\alpha_i \alpha_{i+1} = 0$, $\alpha_i \in \{ 0, 1 \}$,
  be such that $g$ does not start with $f_3 f_5$ nor $f_4 f_6$.

  At first we show that for any $m \ge 6$ and $k \le m-2$, with $(
  k, m ) \neq ( 4, 6 )$, there exists $w = w(k,m)\in F$ such that $w
  f_k f_m = f_m$. If $k = 3$ and $m = 6$ than $w = f_4 f_3 f_1$.  If
  $m = 7$ then all cases of $k$ are covered by the equalities
  \begin{align*}
    f_5 f_4 f_1 \cdot f_5 f_7 &= f_5 f_4 f_2 \cdot f_4 f_7
    = f_5 f_4\cdot f_3 f_7\\
    &= f_5 f_4 f_3^2 f_4 f_6
    = f_5f_4^2f_6 = f_1f_5f_4f_6 = f_1^2f_5f_4f_5 = f_7.
  \end{align*}

  Using Lemma~\ref{lem:reducer}, the element $w(k,m)$ for $m \ge 8$
  may be constructed as
  \[
  w ( k, m ) = f_{m-2} f_{m-3} f_{m-5} f_{m-7} \cdot f_{m-3}
  f_{m-4} \cdots f_k f_{k-1}.
  \]
  Indeed,
  \begin{align*}
    w ( k, m ) \cdot f_k f_m &= f_{m-2} f_{m-3} f_{m-5} f_{m-7}
    \cdot f_{m-3} f_{m-4} \cdots f_k f_{k-1} \cdot f_kf_m\\
    &= f_{m-2} f_{m-3} f_{m-5} f_{m-7} \cdot f_{m-3} f_{m-4}
    \cdots f_k \cdot f_{k+1} f_m \\
    = \cdots &= f_{m-2} f_{m-3} f_{m-5} f_{m-7} \cdot f_{m-2}
    f_m = f_m.
  \end{align*}
  Let $3 \le i_1 < i_2 < \cdots < i_k < n-1$ be the indices such
  that $\alpha_{i_j} = 1$, with the other $\alpha_j = 0$. Let us
  consider the element
  \[
  g_1 = w (i_k, n) w (i_{k-1}, i_k ) \dots w (i_2, i_3) w (i_1,
  i_2).
  \]
  It follows from Lemma~\ref{lem:reducer} that
  \begin{equation}\label{eq:reducer_left}
    g_1 \cdot g = w (i_k, n) w(i_{k-1}, i_k ) \dots w(i_2, i_3) w(i_1,
    i_2) \cdot f_{3}^{\alpha_3} \dots f_{n-2}^{\alpha_{n-2}} f_n = f_n. 
  \end{equation} 

  Now we construct the elements $\ell(g)$ and $r(g)$.  Choose
  \[
  g = f_1^{\alpha_1} f_{3}^{\alpha_3} \cdots f_{n-2}^{\alpha_{n-2}}
  f_n f_{n-1}^{\beta_{n-1}} \cdots f_{2}^{\beta_2} f_1^{\beta_1}.
  \]
  If $g$ includes $z_{n-2}$ then it has normal form $f_1^\alpha
  z_{n-2} f_{n-1}^{\beta_{n-1}} \cdots f_{2}^{\beta_2} f_1^{\beta_1}$,
  and we set $\ell(g) = f_{3-n\%2} f_1^\alpha$ and $r(g) = g_2
  f_{n-1}^{\beta_{n-1}}$. According to the arguments above, we have
  \[
  \ell(g) \cdot g \cdot r(g) = f_{n+1} f_{n-1}^{2 \beta_{n-1}}
  = f_{n+1}.
  \]
  Otherwise $g$ does not include $z_{n-2}$. Assume first that $g$
  starts neither with $f_3 f_5$ nor with $f_4 f_6$. We set $\ell(g) =
  g_1 f_1^{\alpha_1}$ and $r(g) = g_2$. It follows
  from~\eqref{eq:reducer_right} and the arguments above that
  \[
  \ell(g) \cdot g \cdot r(g) = f_n f_{n-1}^{\beta_{n-1}}.
  \]
  If $\beta_{n-1} = 0$, then we are done. If $\beta_{n-1}=1$ it
  follows from Lemma~\ref{lem:idzn} that it suffices to take $\ell(g)
  = fsfs \cdot g_1 f_1^{\alpha_1}$ if $n$ is even, and $\ell(g) =
  sffs \cdot g_1 f_1^{\alpha_1}$ if $n$ is odd; then $\ell(g) \cdot
  g\cdot r(g) = f_{n+1}$.

  The last case is when $g$ starts with $f_1^\alpha f_3 f_5 \cdots
  f_m$ or $f_1^\alpha f_4 f_6 \cdots f_m$, with $6 \le m < n-2$.  We
  may consider respectively the elements $f_2 f_1^\alpha g$ or $f_3
  f_1^\alpha g$, which satisfy the conditions of the previous
  paragraph.
\end{proof}

\end{document}